# A SEMI-ANALYTICAL ANALYSIS OF A FREE CONVECTION BOUNDARY-LAYER FLOW OVER A VERTICAL PLATE

## Haldun Alpaslan PEKER[*] and Galip OTURANÇ

*Department of Mathematics, Faculty of Science, Selcuk University, 42075, Konya, Turkey*

**Abstract.** The aim of this study is to present a reliable combination of the differential transformation method (DTM) and Padé approximants to make, for the first time, a semi-analytic analysis of the problem of free convection boundary-layer flow over a vertical plate as produced by a body force about a flat plate in the direction of the generating body force. Boundary conditions in an unbounded domain, i.e. boundary condition at infinity, pose a problem in general for the numerical solution methods. Without using Padé approximation, the semi-analytical solution obtained by the DTM cannot satisfy infinity boundary conditions. The obtained results are in good agreement with those provided previously by the iterative numerical method. As a result, without taking or estimating missing boundary conditions, the DTM-Padé method provides a simple, non-iterative and effective way for determining the solutions of nonlinear free convection boundary layer problems possessing the boundary conditions at infinity.

**Keywords:** Nonlinear differential equations; Free convection boundary-layer flow; Heat transfer; DTM-Padé method

*AMS classification:* 34A34, 76R10, 76N20, 76N05, 80A20, 34B15, 41A21, 74S30

| Nomenclature | | | |
|---|---|---|---|
| $g$ | Acceleration of gravity | $\theta$ | non-dimensional temperature function |
| $T$ | Temperature | $\alpha$ | thermal diffusivity |
| $T_\infty$ | free stream temperature | $\beta$ | thermal expansion coefficient of fluid |
| $u$ | velocity component along $x$ | $\eta$ | similarity variable |
| $v$ | velocity component along $y$ | $\nu$ | kinematic viscosity |
| $f$ | non-dimensional stream function | | |

---

[*] Corresponding author. E-mail: pekera@hotmail.com

# 1. Introduction

A great number of the nonlinear phenomena can be modeled by nonlinear differential equations in many areas of scientific fields such as engineering, physics, biology, fluid mechanics. Since convection problems come across both in nature and engineering applications, they have attracted a great deal of attention from researches. Free convection flow problems, which results from the action of body forces on the fluid, are one of the common area of interest in the field of convection problems. For the sake of simplicity, the many of free convection boundary-layer problems are considered as a special case of free convection flow about a flat plate parallel to the direction of the generating body force. The most notable model of this topic is the experimental and theoretical considerations of Schmidt and Beckmann [1] concerning the free convection flow of air subject to the gravitational force about a vertical flat plate [2].

This types of free convection boundary-layer problem was studied at the NACA Lewis Laboratory during 1951 and then was analyzed by an iterative numerical method by Ostrach in 1953. After that Na and Habib [3] solved these problems by a parameter method in 1974. In 2005, Kuo [4] employed these problems by the differential transformation method (DTM). In [3] and [4], the boundary conditions of $f''(0)$ and $\theta'(0)$ is taken from the Ostrach [2].

The nonlinear differential equations which have boundary conditions in unbounded domains have a great interest. However, many of the modeled nonlinear equations do not have an analytical solution. Both analytical solutions methods and numerical solutions methods are used to solve these equations. In this study, the DTM is one of the effective and reliable numerical solution method for handling both linear and nonlinear differential equations. In order to overcome the difficulty in unbounded domains, i.e. infinity boundary conditions, the Padé approximant is widely used. In a study [5], the current author and et al. employed combination of the DTM and Padé approximant for a form of Blasius equation. Without using Padé approximant, the analytical solution obtained by DTM cannot satisfy the boundary conditions at infinity.

The present study has an importance due to the fact that this problem is a heat transfer problem consisting boundary conditions at infinity and is solved by a semi-analytical method. In addition, instead of discrete solutions, continuous solutions are obtained without taking initial conditions and/or estimations for the lack of boundary conditions.

## 2. Mathematical Analysis

We consider laminar free-convection flow of an incompressible viscous fluid about a flat plate parallel to the direction of the generating body force. The physical model is shown in Fig.1. By the assumption that the flow in the laminar boundary layer is two dimensional, the continuity equation, the momentum equation, the energy equation and the boundary conditions can be expressed as [1]:

$$\frac{\partial u}{\partial x} + \frac{\partial u}{\partial y} = 0, \qquad (1)$$

$$u\frac{\partial u}{\partial x} + v\frac{\partial u}{\partial y} = g\beta(T - T_\infty) + \nu\frac{\partial^2 u}{\partial y^2}, \qquad (2)$$

$$u\frac{\partial T}{\partial x} + v\frac{\partial T}{\partial y} = \alpha\frac{\partial^2 T}{\partial y^2} \qquad (3)$$

at $x = 0$, $\quad u = 0$ and $T = T_\infty$ $\qquad (4)$

at $y = 0$, $\quad u = v = 0$ and $T = T_0$ $\qquad (5)$

as $y \to \infty$, $\quad u = 0$ and $T = T_\infty$ $\qquad (6)$

where $u$ and $v$ are the velocity in the $x$ and $y$ direction respectively, $\nu$ is the viscosity of the fluid, $\alpha$ is the thermal diffusivity of the fluid.

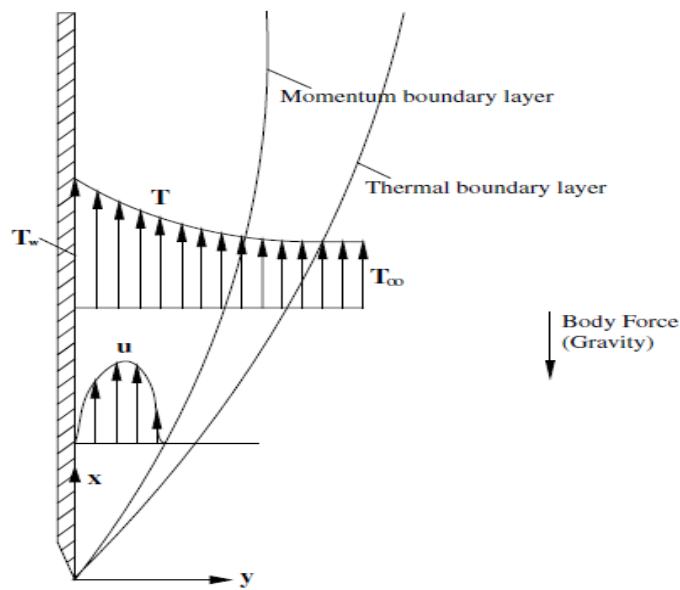

**Figure 1.** Free convection flow over a vertical plate

After a group of transformations [2], the equations (1-3) with respect to the boundary conditions (4-6) reduce to the following form:

$$f'''(\eta) + 3f(\eta)f''(\eta) - 2[f'(\eta)]^2 + \theta = 0 \tag{7}$$

$$\theta''(\eta) + 3\Pr f(\eta)\theta'(\eta) = 0 . \tag{8}$$

The boundary conditions for the equations (7-8) are as follows

$$\begin{array}{llll} at & \eta = 0 & f = f' = 0 \;\; ; & \theta = 1 \\ as & \eta \to \infty & f' = 0 \;\; ; & \theta = 0 \end{array} \tag{9}$$

## 3. Differential Transformation Method

Differential transformation method (DTM) is a semi-analytical method based on Taylor expansion. This method tries to find coefficients of series expansion of unknown function by using the initial data on the problem. By DTM, algebraic equations of initial value and boundary value problems are obtained and then transformation coefficients are calculated. Thus, either closed form series solutions, i.e. analytical solutions or approximate solutions are obtained. This method was firstly introduced in 1954 by Pukhov who especially made some studies related on modelling of electrical and electronical problems. It was applied to electric circuit analysis problems by Zhou [6]. Afterwards, it was applied to several systems and differential equations [8,9,10].

**Definition 1.** The one–dimensional differential transform of a function $y(x)$ at the point $x = x_0$ is defined as follows [7]:

$$Y(k) = \frac{1}{k!}\left[\frac{d^k}{dt^k} y(x)\right]_{x=x_0}, \tag{10}$$

where $y(x)$ is the original function and $Y(k)$ is the transformed function for $k=0,1,2,3,\ldots$ .

**Definition 2.** The differential inverse transform of $Y(k)$ is defined as follows [7]:

$$y(x) = \sum_{k=0}^{\infty} Y(k)(x - x_0)^k . \tag{11}$$

From (10) and (11), we obtain

$$y(x) = \sum_{k=0}^{\infty} \frac{1}{k!} \left[ \frac{d^k}{dt^k} y(x) \right]_{x=x_0} (x-x_0)^k . \qquad (12)$$

Fundamental operations of DTM which can be deduced from Def.1. and Def.2. are given in Table 1.

| Original Function | Transformed Function |
|---|---|
| $y(x) = u(x) \pm v(x)$ | $Y(k) = U(k) \pm V(k)$ |
| $y(x) = cu(x)$ | $Y(k) = cU(k)$ |
| $y(x) = u(x)v(x)$ | $Y(k) = \sum_{r=0}^{k} U(r)V(k-r)$ |
| $y(x) = \dfrac{d^n u(x)}{dx^n}$ | $Y(k) = \dfrac{(k+n)!}{k!} U(k+n)$ |
| $y(x) = u(x)\dfrac{d^2 v(x)}{dx}$ | $Y(k) = U(k) \otimes V(k) = \sum_{r=0}^{k}(k-r+1)(k-r+2)U(r)V(k-r+2)$ |

**Table 1.** The fundamental operations of DTM

## 4. Padé Approximant

A well known fact is that polynomials are used to approximate truncated power series. Further, the singularities of polynomials cannot be seen obviously in a finite plane. Since the radius of convergence of the power series may not large enough to contain the two boundaries, it is not always useful to use the power series [11]. Padé approximants are applied to manipulate the obtained series for numerical approximations to overcome this difficulty. Padé approximant is the best approximation of a polynomial approximation of a function into rational functions of polynomials of given order [12].

Padé approximants are used widely in computer calculations due to the fact that a Padé approximant often gives better approximation of the function than truncating its power series and it may still work where the power series does not converge [12]. Therefore, Padé approximants can be easily computed by using symbolic programming languages such as Maple or Mathematica.

A Padé approximation to *f(x)* on [*a*, *b*] is the quotient of two polynomials, say $P_n(x)$ and $Q_m(x)$ of degrees *n* and *m* respectively, $n, m \in \mathbb{Z}^+$ [12]. The notation $[n/m]$ will be used to denote this quotient. Then the Pade approximant of order $[n/m]$, denoted by $[m,n]_f(x)$ is the following rational function [13]:

$$f(x) \cong R(x) = \frac{Q_m(x)}{P_n(x)} = \frac{\sum_{i=0}^{m} a_i x^i}{1 + \sum_{j=1}^{n} b_j x^j} = \frac{a_0 + a_1 x + a_2 x^2 + \ldots + a_m x^m}{1 + b_1 x + b_2 x^2 + \ldots + b_n x^n}.$$

Note that

$$f(0) = R(0),$$
$$f'(0) = R'(0),$$
$$f''(0) = R''(0),$$
$$\ldots\ldots\ldots$$
$$f^{(m+n)}(0) = R^{(m+n)}(0).$$

In order to obtain better numerical results, the combination of the differential transform method and the diagonal approximants [n/n] will be used.

## 5. Application

Consider a special case of free convection flow about a flat plate parallel to the direction of the generating body force problem given by the equations (7-8), i.e.

$$f'''(\eta) + 3f(\eta)f''(\eta) - 2[f'(\eta)]^2 + \theta = 0$$

$$\theta''(\eta) + 3\Pr f(\eta)\theta'(\eta) = 0.$$

The boundary conditions for the equations (7-8) are as follows

$$\begin{array}{llll} at & \eta = 0 & f = f' = 0 \;;\; & \theta = 1 \\ as & \eta \to \infty & f' = 0 \;;\; & \theta = 0 \end{array}$$

In order to find numerical values of $f''(0)$ and $\theta'(0)$, the series obtained by the differential transform method and the diagonal Padé approximants will be combined.

In this work, we will consider that $Pr = 1$ and also consider that $f''(0) = A$ and $\theta'(0) = B$, where $A$ and $B$ are positive constants.

Now, the differential transform method will be applied to the equations (7) and (8) as follows;

$$(k+1)(k+2)(k+3)F(k+3) = \left(2\sum_{r=0}^{k}(r+1)(k-r+1)F(r+1)F(k-r+1)\right) \\ -\theta(k) - \left(3\sum_{r=0}^{k}\frac{(k-r+1)(k-r+2)F(r)F(k-r+2)}{r!}\right) \quad (13)$$

$$(k+2)(k+1)\theta(k+2) = -3\sum_{r=0}^{k}(k-r+1)F(r)\theta(k-r+1) \ . \quad (14)$$

It is also necessary to find differential transform equivalence form of initial and boundary values. For the initial and boundary values, differential transform equivalences are as follows;

$$F(0) = 0, \ F(1) = 1 \ \text{and} \ F(2) = \frac{A}{2} \quad (15)$$

$$\theta(0) = 0, \ \theta(1) = B \ . \quad (16)$$

Using Maple 15, the following solution is obtained;

$$f(\eta) = \frac{1}{2}A\eta^2 - \frac{1}{6}\eta^3 - \frac{1}{24}B\eta^4 + \frac{1}{48}A^2\eta^5 - \frac{7}{720}A\eta^6 \quad (17)$$

$$\theta(\eta) = 1 + B\eta - \frac{1}{8}AB\eta^4 + \frac{1}{40}B\eta^5 + \frac{1}{240}B^2\eta^6 \ . \quad (18)$$

Now, our aim is to determine a numerical value for $A$ and $B$ using the boundary conditions

$$\lim_{\eta \to -\infty} f'(\eta) = 0 \ \text{and} \ \lim_{\eta \to -\infty} \theta(\eta) = 0 \quad (19)$$

which is in condition (9). In order to do this, derivative of polynomial solution (17) should be taken. After that the formation of this equation and the equation (18) by Padé approximation, we will apply the condition (19) to the obtained rational functions.

The [3/3] diagonal Padé approximant is computed by using Maple 15 for which the maple codes are shown in Table 2.

```
restart:
Digits:=10:
m:=6: #Order
F[0]:=0:
F[1]:=0:
F[2]:=A/2:
T[0]:=1:
T[1]:=B:
for k from 0 to m+2 do
F[k+3]:=(2*sum((r+1)*F[r+1]*(k-r+1)*F[k-r+1],r=0..k)-T[k]-
3*sum((k-r+1)*(k-r+2)*F[r]*F[k-r+2]/r!,r=0..k))*k!/(k+3)!:
T[k+2]:=-3*sum((k-r+1)*F[r]*T[k-r+1],r=0..k)*k!/(k+2)!:
od:
f:=0:
t:=0:
for k from 0 to m do
f:=f+F[k]*x^k:
t:=t+T[k]*x^k:
od:
print(f):
print(t):
with(numapprox):
pade(t,x,[3,3]);
pade(diff(f,x),x,[3,3]);
solve({limit(pade(diff(f,x),x,[3,3]),x=infinity)=0.,limit(pade
(t,x,[3,3]),x=infinity)=0.},[A,B]);
```

**Table 2.** Maple codes

According to the above procedure, obtained values for *A* and B are listed in Table 3 with the ones obtained in [2].

| Padé approximant | A | A (Ref. [2]) | B | B (Ref.[2]) |
|---|---|---|---|---|
| [3/3] | 0.5506447081 | 0.6421 | -0.8654409691 | -0.5671 |

**Table 3.** Padé approximants and numerical values of *A* and *B*

In addition to these, there does not exist analytical solution for this problem. Therefore, in order to see the consistency of solutions, the series solutions up to the sixth term of the problem are calculated with respect to the numerical values of *A* and *B* for which $A = 0.5506447081$ and $B = -0.8654409691$ according to the [3/3] diagonal Padé approximant. Furthermore, they are calculated with step size 0.1 on the interval [0,1] and the results are shown in Fig.2 and Fig.3.

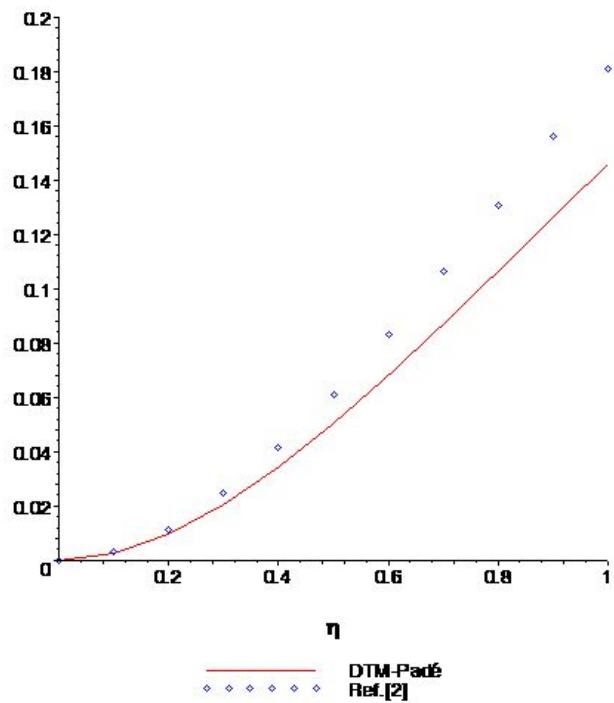

**Figure 2.** DTM-Padé and Ref. [2] solutions of $f$ for the numerical values of *A* and *B* for the [3/3] diagonal Padé approximant for *Pr*=1

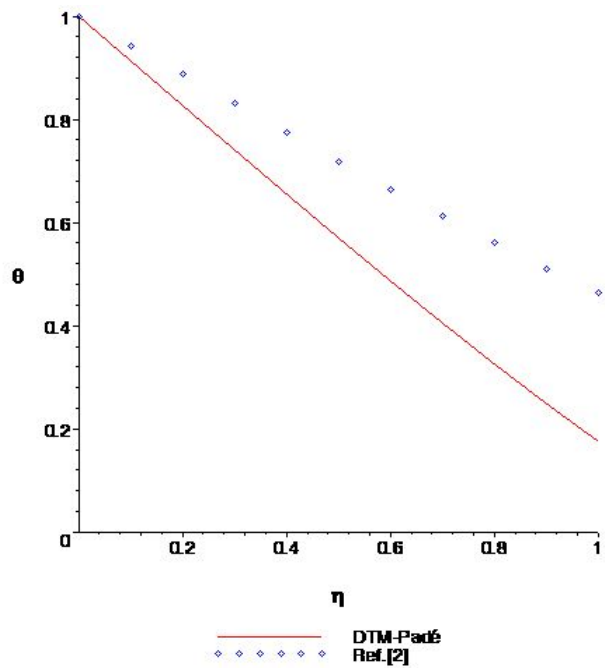

**Figure 3.** DTM-Padé and Ref. [2] solutions of $\theta$ for the numerical values of *A* and *B* for the [3/3] diagonal Padé approximant for *Pr*=1

## 6. Conclusion

Boundary conditions at infinity pose a problem in some of heat transfer problems as well as many other nonlinear differential equations problems. Analytical solution of the problem cannot be obtained due to both the lack of linearity of the problem and the boundary conditions at infinity. In two-point boundary value problems where one point is infinity, two questions arise. The first one is that where infinity is and the second question is that when a satisfactory approximation to a solution has been obtained. In this study, a special case of free convection flow about a flat plate parallel to the direction of the generating body force problem, which is a nonlinear problem and possessing the boundary conditions at infinity, is solved semi-analytically for the first time without taking or estimating boundary conditions $f''(0)$ and $\theta'(0)$. It is easily seen that the solutions in the literature were made by using the boundary conditions for $f''(0)$ and $\theta'(0)$ given by Ostrach [2]. Even Ostrach himself obtained his solutions by estimating the values for $f''(0)$ and $\theta'(0)$ [2]. In addition to these, it is clear that the obtained series is convergent.

The obtained results by the DTM-Padé are compared with, for the first time, the ones found in [2].

**Acknowledgement**

This research is supported by TUBITAK (The Scientific and Technological Research Council of Turkey) and Selcuk University Scientific Research Project Coordinatorship (BAP). This study is a part of the corresponding author's Ph.D. Thesis.